\documentclass{amsart}

\usepackage{amssymb,amsmath, amsthm, amsfonts}

\newcommand{\rd}{\mathbb{R}^d}
\newcommand{\rtwo}{\mathbb{R}^2}
\newcommand{\rthree}{\mathbb{R}^3}
\newcommand{\R}{\mathbb{R}}

\newcommand{\zed}{\mathbb{Z}}
\newcommand{\pq}{L^p\to L^q}

\newcommand{\ca}{\mathcal{A}}

\newtheorem{thm}{Theorem}
\newtheorem{cor}{Corollary}
\newtheorem{lem}{Lemma}

\begin{document}

\title[$L^p$ improving along polynomial curves]{Universal $L^p$ improving for averages along\\ polynomial curves in low dimensions}

\author{Spyridon Dendrinos, Norberto Laghi and James Wright}

\address{\flushleft Dendrinos:
\newline
Department of Mathematics, University of Bristol, University Walk,
Bristol, BS8 1TW, United Kingdom} \email{S.Dendrinos@bris.ac.uk}
\address{\flushleft Laghi and Wright:
\newline
Maxwell Institute for Mathematical Sciences, The University of
Edinburgh, JCMB, King's Buildings, Edinburgh EH9 3JZ, United
Kingdom} \email{N.Laghi@ed.ac.uk} \email{J.R.Wright@ed.ac.uk}


\thanks{N. Laghi and J. Wright were supported in part by an EPSRC grant}

\begin{abstract}
We prove sharp $\pq$ estimates for averaging operators along
general polynomial curves in two and three dimensions. These
operators are translation-invariant, given by convolution with the
so-called affine arclength measure of the curve and we obtain
universal bounds over the class of curves given by polynomials of
bounded degree. Our method relies on a geometric inequality for
general vector polynomials together with a combinatorial
argument due to M. Christ. Almost sharp Lorentz space estimates
are obtained as well.
\end{abstract}

\maketitle

\section{Introduction and statement of results}

Recently there has been considerable attention given to certain
euclidean harmonic analysis problems associated to a curve or
surface where the underlying euclidean arclength or surface measure (which
typically defines the classical problem) is replaced by the
so-called affine arclength or surface measure. This has the effect
of making the problem affine invariant as well as invariant under
reparametrisations of the underlying variety. For this reason
there have been many attempts to obtain universal results,
establishing uniform bounds over a large class of curves or
surfaces. The affine arclength or surface measure also has the
mitigating effect of dampening any curvature degeneracies of the
curve or surface and therefore the expectation is that the
universal bounds one seeks will be the same as those arising from
the most non-degenerate situation.

This line of research has been actively pursued for the
problem of Fourier restriction, a central problem in euclidean
harmonic analysis; see for example \cite{bassam-1}, \cite{bos}, \cite{ckz},
\cite{cz}, \cite{DW}, \cite{drury-90}, \cite{dm-1}, \cite{dm-2},
\cite{oberlin-sjolin}, \cite{oberlin-surfaces} and
\cite{bassam-2}. Drury initiated an investigation along
these lines for the problem
of achieving precise regularity
results for averages along curves or surfaces, in particular
determining sharp $\pq$ estimates, and this has been followed
up by several authors; see for example \cite{1}, \cite{choi-03},
\cite{drury-90}, \cite{oberlin-91},
\cite{oberlin-smoothing}, \cite{oberlin-polynomial}, \cite{12},
\cite{pan-93},
\cite{pan-94} and \cite{pan-96}.

In this paper we continue an investigation by Oberlin
to establish such a result for averaging operators along general
polynomial curves in ${\mathbb R}^d$ when $d=2$ or $d=3$ (in
\cite{oberlin-polynomial}, the $d=2$ case was fully resolved and
partially resolved for $d=3$). More specifically, if $\gamma: I
\rightarrow {\mathbb R}^d$ parametrises a smooth curve in
${\mathbb R}^d$ on an interval $I$, set
$$
L_{\gamma}(t) = {\rm det} (\gamma^{\prime}(t) \cdots
\gamma^{(d)}(t) );
$$
this is the determinant of a $d\times d$ matrix whose $j$th column
is given by the $j$th derivative of $\gamma$, $\gamma^{(j)}(t)$.
The affine arclength measure $\nu = \nu_{\gamma}$ on $\gamma$ is
defined on a test function $\phi$ by
$$
\nu (\phi) = \int_I \phi(\gamma(t))
|L_{\gamma}(t)|^{\frac{2}{d(d+1)}} dt;
$$
one easily checks that this measure is invariant under
reparametrisations of $\gamma$. A basic problem in the theory of
averaging operators along curves (or more generally, for
generalised Radon transforms) is to determine the exponents $p$
and $q$ so that the apriori estimate
\begin{equation}\label{apriori} \|Tf \|_{L^q({\mathbb R}^d)} \le C
\|f\|_{L^p({\mathbb R}^d)}
\end{equation}
holds uniformly for a large class of curves $\gamma$ where
$$
Tf(x)=f*\nu(x)=\int_{I}f(x-\gamma(t))|L_{\gamma}(t)|^{\frac{2}{d(d+1)}}dt.
$$

Use of the affine arclength measure allows us to think about
global estimates, not only establishing \eqref{apriori} with a
constant $C$ uniform over a large class of curves but also
possibly obtaining such a constant independent of the
parametrising interval $I$. On the other hand thinking of $T$ as a
local operator and thus insisting $T$ preserve all $L^p$ spaces,
the constant $C$ in \eqref{apriori} will then necessarily depend
on $I$. As discussed above, the exponents $p$ and $q$ in
\eqref{apriori} that we expect should come from the most
non-degenerate situation which in this case is the curve
$\gamma(t) = (t,\ldots,t^d)$ in ${\mathbb R}^d$ where $L_{\gamma}
\equiv constant$. With regards to local estimates in this case
(and thus allowing $C$ to depend on $I$), by testing
\eqref{apriori} on $f = \chi_{B_{\delta}}$ where $B_{\delta}$ is
the ball of radius $\delta$ with centre 0, $f = \chi_{D_{\delta}}$
where $D_{\delta} = \{|x_1|\le \delta, \ldots, |x_d| \le \delta^d
\}$ and using duality, one easily sees that the exponents $p$ and
$q$ necessarily satisfy
$$
(1/p,1/q)\in H_d=\text{hull}\{(0,0),(1,1),A_d,B_d\}, \ \ {\rm
where} \ A_d=\bigl(2/(d+1),(2d-2)/(d^2+d)\bigr)
$$
and $B_d=\bigl((d^2-d+2)/(d^2+d),(d-1)/(d+1)\bigr)$. It is a
remarkable result of Christ \cite{2} that (up to the endpoints
$A_d$ and $B_d$) these restrictions on $p$ and $q$ are in fact
sufficient for \eqref{apriori} to hold in this non-degenerate
situation. It is our understanding that Stovall \cite{13},
building on an argument of Christ \cite{3}, has converted Christ's
restricted weak-type estimates at $A_d$ and $B_d$ into strong type
estimates. With regards to global estimates in this non-degenerate
situation $\gamma(t) = (t,\ldots, t^d)$ (ensuring $C$ in
\eqref{apriori} can be taken to be independent of $I$), by a simple
scaling argument or by taking $f = \chi_{D_{\delta}}$ but now
letting $\delta$ vary over all the positive reals, one sees that
necessarily we must have $1/q = 1/p - 2/d(d+1)$. Furthermore, the
necessary conditions for the local estimates give us the added
restriction $(d^2 + d)/(d^2 - d +2) \le p \le (d+1)/2$.

To date, progress that has been made to establish universal bounds
in \eqref{apriori} for curves $\gamma$ where $L_{\gamma}
\not\equiv constant$ has not been as substantial as for the
corresponding problem of Fourier restriction. The case for curves
$\gamma(t) = (t,\phi(t))$ given as the graph of a convex function
$\phi$ has been considered by Choi, Drury, Oberlin and Pan and the best result
here is due to Oberlin \cite{oberlin-smoothing} where the
additional hypothesis that $\phi''$ is monotone increasing is
imposed and then only a weak-type estimate is obtained
at the endpoint $(2/3, 1/3)$ (in \cite{1} Choi obtained strong
type estimates at $(2/3,1/3)$ but these estimates are not
universal -- the constant $C$ in \eqref{apriori} depends on $\phi$
-- and in fact the author needs to impose much more stringent
conditions on $\phi$).

Compare this with the situation for the corresponding Fourier
restriction problem in two dimensions where Sj\"olin \cite{sjolin}
obtained uniform bounds over the class of {\it all} convex curves
-- see also \cite{oberlin-sjolin}. The class of convex curves is a
natural class to examine in light of simple counterexamples to
\eqref{apriori} where $L_{\gamma}$ changes sign too often (of
course if $\gamma$ is convex, $L_{\gamma}$ does not change sign).
By the above discussion on necessary conditions, we see that the
endpoint estimate to aim for in \eqref{apriori} is $(2/3,1/3)$ in
two dimensions. Consider the curve $\gamma$ given by $\gamma(t) =
(t, t^k \sin(1/t) )$. By testing \eqref{apriori} on $f =
\chi_{D_{\delta}}$ where $D_{\delta} = \{(x,y) : |x|\le \delta,
|y| \le \delta^k \}$ one easily shows that if \eqref{apriori} were
to hold for this example, then $1/q \ge 1/p - (k-1)/3(k+1)$.
Therefore if $L_{\gamma}$ changes sign too often then
\eqref{apriori} may not hold uniformly for all curves in the
expected $L^p$ range.

In \cite{oberlin-polynomial} Oberlin established \eqref{apriori}
in two dimensions for the family of polynomial curves $\gamma(t) =
{\bf P}(t) = (P_1(t), P_2(t))$ where each $P_1$ and $P_2$ is a
general real polynomial of bounded degree.
Specifically he established \eqref{apriori} with a constant $C$
only depending on the the degrees of the polynomials defining
${\bf P}$. This is a natural class of curves to consider as the
number of sign changes of $L_{{\bf P}}$ is controlled by the
degree of the polynomials $P_j$. Furthermore Oberlin established
\eqref{apriori} in three dimensions for polynomial curves of the
form ${\bf P}(t) = (t, P_2(t), P_3(t))$ but the estimates are not
universal in the sense that the constant $C$ can be taken to
depend only on the degrees of the polynomials. For the
corresponding Fourier restriction problem in the setting of
polynomial curves, see \cite{bos} and \cite{DW}.

In this paper we give an alternative approach to the results in
\cite{oberlin-polynomial} and strengthen the three dimensional
result to general polynomial curves ${\bf P}(t) = (P_1(t), P_2(t),
P_3(t))$; furthermore all estimates will be uniform over the class
of polynomials of bounded degree. Our hope is that this approach
will generalise to general polynomials curves in all
dimensions.

From now on we shall focus on the operator
\begin{equation}\label{ca}\ca f(x)=\int_I
  f(x-\mathbf{P}(t))\left|L_{\mathbf{P}}(t)\right|^{\frac{2}{d(d+1)}}dt.
\end{equation}
We are now ready to state our main result which is a global
estimate.

\begin{thm}\label{main} Let $d=2,3.$ Then for every $\epsilon>0,$
\[
\|\ca f\|_{L^{\frac{d^2+d}{2d-2},\frac{d+1}{2}+\epsilon}(\rd)}\leq
C \|f\|_{L^{\frac{d+1}{2}}(\rd)}\] and \[\|\ca
f\|_{L^{\frac{d+1}{d-1},\frac{d^2+d}{d^2-d+2}+\epsilon}(\rd)}\leq
C \|f\|_{L^{\frac{d^2+d}{d^2-d+2}}(\rd)},
\]
where the constant $C$ depends only on $\epsilon>0$, the degrees of
the polynomials defining the curve $\mathbf{P}$ and in particular
not on the parametrising interval $I$.
\end{thm}
When $d=2$ there is just a single endpoint and the above two
estimates agree. Here $L^{p,r}(\rd)$ denote the familiar Lorentz
spaces. Since $C$ can be taken to be independent of $I$ and $\ca$
is a positive operator, Theorem \ref{main} is equivalent to
establishing the concluding estimates for the global analogue of
$\ca$ where the integration in \eqref{ca} is replaced by the
entire real line.

Utilising Theorem \ref{main} and the well-known local estimates
giving boundedness for our operators on the line $p=q,$ we obtain
the following consequence.
\begin{cor}\label{cor-main} Let $d=2,3.$ Then if $(1/p,1/q)\in H_d$,
\[\|\ca f\|_{L^p(\rd)}\leq C\|f\|_{L^q(\rd)},\]
where the bound $C$ depends only on the degrees of the polynomials
defining the curve $\mathbf{P}$ and on the interval $I.$
\end{cor}

The proof of Theorem \ref{main} combines an elegant
combinatorial argument of Christ in \cite{2}, together with a
recent geometric inequality for vector polynomials which was
established in \cite{DW}. Christ's method is elementary but powerful
and has seen applications outside the model curve
case $(t,\ldots, t^d)$ (see \cite{bcw}, \cite{4} and \cite{5}) as
well as substantial generalisations (see \cite{CE} and \cite{14}).
We mention again that Christ has developed a method that may be
used to deduce strong-type estimates (even Lorentz type estimates)
from restricted weak-type estimates (see \cite{3}) and we will
follow this method to deduce the Lorentz bounds in Theorem
\ref{main}.

Finally, we wish to emphasise the fact that the result of Theorem
\ref{main} is obtained by using slightly different ingredients in
different dimensions; whilst the basic techniques employed do not
change, the relevant arguments need to be suitably adjusted. This
is reflected in the structure of the paper: in the next section we
recall the rudiments of Christ's argument in \cite{2} followed by
a description in \S 3 of the {\it key} geometric inequality for
polynomial curves established in \cite{DW}, an essential fact in
our arguments. In \S 4 we deal with the restricted weak-type
estimates in three dimensions, and in \S 5 we show how these
can be turned into strong-type and indeed Lorentz-space estimates,
again in three dimensions. In \S 6 we produce the necessary
arguments needed to deal with the two-dimensional case, while in
the last section we shall discuss the sharpness of our main
result.

\emph{Notation.} Throughout this paper, whenever we write
$A\lesssim B$ or $A = O(B)$ for any two nonnegative quantities $A$
and $B,$ we mean that there exists a strictly positive
constant $c,$ possibly depending on the degree of the map
$\mathbf{P},$ so that $A\leq cB;$ this constant is subject to
change from line to line and even from step to step. We also write
$A\sim B$ if $A \lesssim B \lesssim A$.

\section{Rudiments of Christ's argument}
For a nonnegative finite measure $\mu$ supported on an interval
$I$ and a curve parametrised by $\gamma : I \to {\mathbb R}^d$,
consider the averaging operator
$$
Af(x) \ = \int f(x-\gamma(t)) \, d\mu(t).
$$
In this section we recall the basics of the combinatorial argument
of Christ in \cite{2} to prove a restricted weak-type estimate $A
: L^{p,1}({\mathbb R}^d)\to L^{q,\infty}({\mathbb R}^d)$. This is
equivalent to proving
\begin{equation}\label{rw}
\langle A\chi_E,\chi_F\rangle \ \lesssim \ |E|^{1/p}|F|^{1/q'}
\end{equation}
for any two measurable sets $E,F \subset {\mathbb R}^d$ where $|\cdot
|$ denotes the Lebesgue measure. Without loss of generality we may
assume that $|E|$, $|F|$ and $\langle A\chi_E,\chi_F\rangle$ are
all positive quantities. Define two positive parameters $\alpha$
and $\beta$ by the relations
$$
\alpha := \frac{1}{|F|}\langle A\chi_E,\chi_F\rangle, \ \ \beta :=
\frac{1}{|E|}\langle A^{*}\chi_F, \chi_E\rangle \ \ \ {\rm so \
that} \ \ \alpha |F| = \beta |E|
$$
where $A^{*}f (y) = \int f (y + \gamma(t)) \, d\mu(t)$. Thus
$\alpha$ is the average value of $A\chi_E$ on $F$ and $\beta$ is
the average of $A^{*}\chi_F$ on $E$.

By passing to refinements of the sets $E$ and $F$, without
changing significantly the basic quantity $ K := \langle
A\chi_E,\chi_F\rangle = \langle \chi_E,A^{*}\chi_F\rangle$ to be
estimated in \eqref{rw}, we will be able to bound pointwise
$A\chi_E$ by $\alpha$ on $F$ and bound pointwise  $A^{*}\chi_F$ by
$\beta$ on $E$. Precisely one defines the following refinements of
$E$ and $F$:
$$
F_1 = \{x\in F: A\chi_E(x) \ge \alpha/2\}, \ \ \ E_1 = \{y\in E:
A^{*}\chi_{F_1}(y) \ge \beta/4 \},
$$
$$
F_2 = \{x\in F_1: A\chi_{E_1}(x) \ge \alpha/8\},  \ldots, \ E_n =
\{y\in E_{n-1}: A^{*}\chi_{F_n}(y)\ge \beta/2^{2n}\},
$$
etc... It is a simple matter to check that $\langle A\chi_{E_n},
\chi_{F_n} \rangle \ge K/2^{2n}$ and $\langle \chi_{E_n},
A^{*}\chi_{F_{n+1}} \rangle \linebreak \ge K/2^{2n+1}$ for each $n$ and so
$E_n, F_n \not= \emptyset$.

If $d=3$, we fix an $x_0 \in F_{2}$, set $S = \{ s \in I: x_0 -
\gamma(s) \in E_{1} \}$ and note
\begin{equation}\label{S-3}
\mu(S) \ = \ A\chi_{E_1}(x_0) \ \ge \ \alpha/8.
\end{equation}
Next observe that for every $s\in S$, if $T_s = \{ t\in I: x_0 -
\gamma(s) + \gamma(t) \in F_1 \}$, then
\begin{equation}\label{T_s-3}
\mu(T_s) \ = \ A^{*}\chi_{F_1}(x_0 - \gamma(s)) \ \ge \ \beta/4.
\end{equation}
Finally we see that for every $s\in S$ and $t\in T_s$, if $U_{s,t}
= \{u\in I: x_0 - \gamma(s) + \gamma(t) - \gamma(u) \in E\}$, then
\begin{equation}\label{U_{s,t}-3}
\mu(U_{s,t}) \ = \ A\chi_E (x_0 - \gamma(s) + \gamma(t)) \ \ge \
\alpha/2.
\end{equation}
Hence we end up with a structured parameter domain ${\mathcal P} =
\{(s,t,u)\in I^3: s\in S, t\in T_s, u\in U_{s,t} \}$ so that if
$\Phi_{\gamma} (s,t,u) := x_0 - \gamma(s) + \gamma(t) - \gamma(u)$,
$\Phi_{\gamma}({\mathcal P}) \subset E$. Therefore if
$\Phi_{\gamma}$ is injective we have
$$
|E| \ \ge \ \mathop {\int\!\!\!\int\!\!\!\int}_{\mathcal P}
|J_{\Phi_{\gamma}}(s,t,u)| ds dt du \ = \ \int_S
\int_{T_s}\int_{U_{s,t}} |J_{\Phi_{\gamma}}(s,t,u)| ds dt du
$$
where $J_{\Phi_{\gamma}}(s,t,u) = {\rm det} ({\gamma}'(s) \
{\gamma}'(t) \ \gamma'(u) )$ is the determinant of the Jacobian
matrix for the mapping $\Phi_{\gamma}$, reducing matters to
understanding the smallness of $J_{\Phi_{\gamma}}$ (for instance,
sublevel sets of $J_{\Phi_{\gamma}}$) in order to bound from below
the above integral over the structured set ${\mathcal P}$. If
$\gamma(t) = (t,t^2,t^3)$ (the non-degenerate example in three
dimensions) and $\mu = |\cdot |$ is Lesbesgue measure, then simply
$J_{\Phi_{\gamma}}(s,t,u) = 6(s-t)(t-u)(s-u)$ and so \eqref{S-3},
\eqref{T_s-3} and \eqref{U_{s,t}-3} quickly imply $|E|\ge
\beta^2 \alpha^4$ which gives \eqref{rw} with $p=2$ and $q=3$, the
desired endpoint estimate in this case.

If $d=2$, we fix a $y_0 \in E_{1}$, set $S = \{ s \in I: y_0 +
\gamma(s) \in F_{1} \}$ and note
\begin{equation}\label{S-2}
\mu(S) \ = \
A^{*}\chi_{E_1}(y_0) \ \ge \ \beta/4.
\end{equation}
Next observe that for every $s\in S$, if $T_s = \{ t\in I: y_0 +
\gamma(s) - \gamma(t) \in E\}$, then
\begin{equation}\label{T_s-2}
\mu(T_s) \ = \ A\chi_{E}(y_0 + \gamma(s)) \ \ge \ \alpha/2.
\end{equation}
Hence we end up with a structured parameter domain ${\mathcal P} =
\{(s,t)\in I^2: s\in S, t\in T_s \}$ so that if $\Phi_{\gamma}
(s,t) := y_0 + \gamma(s) - \gamma(t)$, $\Phi_{\gamma}({\mathcal
P}) \subset E$. Therefore if $\Phi_{\gamma}$ is injective we have
$$
|E| \ \ge \ \mathop {\int\!\!\!\int}_{\mathcal P}
|J_{\Phi_{\gamma}}(s,t)| \, ds\, dt \ = \ \int_S \int_{T_s}
|J_{\Phi_{\gamma}}(s,t)| \, ds\, dt
$$
where $J_{\Phi_{\gamma}}(s,t) = - {\rm det} ({\gamma}'(s) \
{\gamma}'(t))$. If $\gamma(t) = (t,t^2)$ (the non-degenerate
example in two dimensions) and $\mu = |\cdot |$ is Lesbesgue
measure, then $J_{\Phi_{\gamma}}(s,t) = 2(s-t)$ and so
\eqref{S-2}, \eqref{T_s-2} imply $|E|\ge \beta \alpha^2$ which
gives \eqref{rw} with $p=3/2$ and $q=3$, the desired endpoint
estimate in this case.

Interestingly when we consider a general polynomial curve
$\gamma(t) = {\bf P}(t) = (P_1(t), \linebreak P_2(t))$ in two dimensions with
$\mu$ the affine arclength measure on ${\bf P}$,
we will only be able to prove
$$
\mathop {\int\!\!\!\int}_{\mathcal P} |J_{\Phi_{\gamma}}(s,t)| \,
ds\, dt \ = \ \int_S \int_{T_s} |J_{\Phi_{\gamma}}(s,t)| \, ds\,
dt \ \ge \ \beta \alpha^2
$$
in the range $\alpha \le \beta$. In fact, without further
information,  this integral bound is false in general in the range
$\beta \le \alpha$. Nevertheless, due to the fact that the sharp
endpoint estimate lies on the line of duality $L^p \to L^{p'}$, it
will be the case that $|E|\ge \beta \alpha^2$ for {\it all}
$\alpha, \beta$. The failure of this integral bound in the range
$\beta \le \alpha$ leads to some further difficulties when
establishing the Lorentz bounds and these difficulties do not
present themselves in the three dimensional case. This is why we
choose to address the three dimensional case first.


\section{A geometric inequality}
As we have seen in the previous section, Christ's argument in
\cite{2} is based in part on analysis of the map
$$
\Phi_{{\bf P}}(t_1,\ldots,t_d) = (-1)^d {\bf P}(t_1) + (-1)^{d+1} {\bf P}(t_2) +
\cdots - {\bf P}(t_d).
$$
In particular it would be desirable to have the following
properties about $\Phi_{\bf P}$:
\vspace{1em}

{\bf Key properties}
\\
\hspace*{.25in} (a) \ $\Phi_{\bf P}$ is 1-1; \\
\hspace*{.25in} (b) \ $|J_{\Phi_{\bf P}}(t_1,\ldots,t_d)| \ge C
\prod_{j=1}^{d} |L_{\bf P}(t_j)|^{\frac{1}{d}} \prod_{j<k} | t_j - t_k |$\\
{\it where $J_{\Phi_{\bf P}}(t_1,\ldots,t_d) = \pm {\rm det} ({\bf
P}'(t_1) \cdots {\bf P}'(t_d))$ is the determinant of the Jacobian
matrix for the mapping $\Phi_{\bf P}$ and $L_{\bf P}(t) = {\rm
det} ({\bf P}'(t) \cdots {\bf P}^{(d)}(t))$ was introduced in the
introduction as part of the definition of the affine arclength
measure along ${\bf P}$. }
\vspace{1em}

As we have seen the injectivity of $\Phi_{\bf P}$ allows us to
reduce matters to examining integrals of $J_{\Phi_{\bf P}}$ over
various structured sets of $(t_1,\ldots,t_d)$. And then the
geometric inequality, property (b), will make the examination of
these integrals feasible. Even in the non-degenerate case ${\bf
P}(t) = (t,t^2,\ldots,t^d)$, $\Phi_{\bf P}$ is not quite 1-1 but
it is $d!$ to 1 off a set of measure zero. Furthermore in this
case, the geometric inequality (b) is an equality.

For polynomial curves both (a) and (b) are false in general.
However in \cite{DW}, a collection of $O(1)$ disjoint open
intervals $\{I\}$ was found which decomposes ${\mathbb R} = \cup
{\overline I}$ so that on each $I^d$, $\Phi_{\bf P}$ is $d!$ to 1
off a set of measure zero and the geometric inequality (b) holds.
With this decomposition we will restrict our original operator
$\ca$ to each interval $I$ and apply Christ's argument. The
decomposition is valid only under the assumption that $L_{\bf P}
\not\equiv 0$. Of course if $L_{\bf P} \equiv 0$, then the
estimates in \eqref{apriori} are trivial and so, without loss of
generality, the non-degeneracy assumption $L_{\bf P} \not\equiv 0$
will be in force for the remainder of the paper.

The decomposition is produced in two stages. The first stage
produces an elementary decomposition of ${\mathbb R} = \cup
{\overline J}$ so that on each open interval $J$, various
polynomial quantities (more precisely, certain determinants of
minors of the $d\times d$ matrix $({\bf P}'(t) \cdots {\bf
P}^{(d)})$, including $L_{\bf P}$) are single-signed. This allows
us to write down a formula relating $J_{\Phi_{\bf P}}$ and $L_{\bf
P}$. When $d=2$ this formula is particularly simple; namely,
$$
J_{\Phi_{\bf P}}(s,t) = P^{\prime}_1(s) P^{\prime}_1(t) \int_s^t
\frac{L_{\bf P}(w)}{P^{\prime}_1(w)^2} dw
$$
for any $s,t \in J$ (here ${\bf P} = (P_1, P_2)$). From this, one
can establish the injectivity of $\Phi_{\bf P}$ on $\{
(t_1,\ldots,t_d) \in J^d : \ t_1 < \cdots < t_d \}$. Next we
decompose each ${\overline J} = \cup {\overline I}$ further so
that on each open interval $I$, (b) holds. More precisely, we have inequality (b)
for all $(t_1,\ldots,t_d)\in I^d$ where $C$ depends only on $d$
and the degrees of the polynomials defining ${\bf P}$. 

This second stage decomposition ${\overline J} = \cup {\overline
I}$ is more technical and derived from a certain algorithm which
uses two further decomposition procedures generated by individual
polynomials. These further decomposition procedures are used in
tandem and have the effect of reducing (\ref{ca}) to open intervals $I$ on
which various polynomials, including $L_{\bf P}$, behave like a
centred monomial. Furthermore the algorithm exploits in a crucial
way the affine invariance of the inequality (b);
that is, the inequality is invariant under replacement of ${\bf P}$ by $A{\bf P}$ for any invertible $d\times d$ matrix $A$.

To recapitulate, in \cite{DW} a decomposition ${\mathbb R} = \cup
{\overline I}$ where $\{I\}$ is an $O(1)$ collection of open
disjoint intervals was produced so that the following three
properties hold for each $I$:

\vspace{1em}

\noindent (P1) the map $\Phi_{\bf P}$ is 1-1 on the region
$D = \{(t_1,\ldots,t_d) \in I^d : \, t_{1} < \cdots < t_{d} \};$
\\
(P2) for $t\in I$, $|L_{\bf P}(t)| \sim A_I |t - b_I|^{k_I}$ for some
$A_I > 0$, $b_I \notin I$ and integer $k_I \ge 0$;
\\
(P3) for $(t_1,\ldots,t_d)\in I^d$,
$$
|J_{\Phi_{\Gamma}}(t_1,\ldots,t_d)| \ge C \prod_{j=1}^{d}
|L_{\Gamma}(t_j)|^{\frac{1}{d}} \prod_{j<k} | t_j - t_k |
$$
where $C$ depends only on $d$ and the degrees of the polynomials
defining ${\bf P}$.

\section{Restricted weak-type estimates}
As mentioned above it suffices to carry out our analysis for the
globally defined operator
\begin{equation}\label{5}\ca_\R f(x)=\int_{\R} f(x-\mathbf{P}(t))
\left|L_{\mathbf{P}}(t)\right|^{\frac{2}{d(d+1)}}dt,\end{equation}
and we begin by proving the desired restricted weak-type
estimates. We have the following.
\begin{thm}\label{restricted} Let $d=3;$ the operator (\ref{5}) satisfies
\begin{align}\label{6} \ca_\R :L^{2,1}(\rthree)\to L^{3,\infty}(\rthree),\\
\label{7}\ca_\R:L^{3/2,1}(\rthree)\to L^{2,\infty}(\rthree),\end{align}
where the bounds depend only on the degree of $\mathbf{P}.$
\end{thm}
\noindent {\bf Proof } By duality it suffices to establish just one of these estimates,
say \eqref{6}, and as we have seen in \S 2, this in turn is
equivalent to proving
\begin{equation}\label{rw-3}
\langle \ca_\R \chi_E,\chi_F\rangle \ \lesssim \ |E|^{1/2}|F|^{2/3}
\end{equation}
for all pairs of measurable sets $E,F\subset {\mathbb R}^3$. We now
apply the decomposition procedure described in \S 3 to the vector
polynomial ${\bf P}(t) = (P_1(t),P_2(t),P_3(t))$, decomposing $\R
= \cup {\overline I}$ into $O(1)$ disjoint open intervals $\{I\}$
so that for each $I$, properties (P1), (P2) and (P3) hold.

Thus for each $I$, we need only consider the operator
\[\ca_{I}f(x)=\int_I f(x-\mathbf{P}(t))|t-b|^{k/6}dt := \int_I
f(x-\mathbf{P}(t))d\mu(t),\] and prove \eqref{rw-3} for $\ca_I$,
uniformly in $I$. Here $b = b_I \notin I$, $k = k_I$ is some
nonnegative integer and $\mu = \mu_I$  \footnote{It will be
helpful, for the calculations that will follow, to keep in mind
that the $\mu$ measure of a measurable set $J\subset I$ is given
by $\int_J |t-b|^{k/6}dt.$} is a measure supported in $I$.
Introducing the positive parameters $\alpha = \alpha_I$ and $\beta
= \beta_I$ as in \S 2, we see that
\begin{equation}\label{rw-E}
\left|\langle\ca_I \chi_E,\chi_F\rangle\right|\lesssim
|E|^{1/2}|F|^{2/3}\iff |E|\gtrsim \alpha^4\beta^2,
\end{equation}
uniformly in $I$. From \S 2, we see that there is a point $x_0\in
F$ and
\begin{align*}& S \subset I\text{ so that
  }\mu(S)\gtrsim\alpha; \\
& \text{for each }s\in S \text{ there is a }T_{s}\subset
  I\text{ so that }\mu(T_{s})\gtrsim\beta;\\
&\text{for each }t\in T_{s}\text{ there is a }U_{s,t}\subset I
\text{ so that } \mu(U_{s,t})\gtrsim\alpha;\\
&\text{if }\mathcal{P}=\left\{(s,t,u)\in I^3: s\in  S,t\in T_{s},
u\in U_{s,t}\right\}\text{ then } x_0+\Phi_{\mathbf{P}}
(\mathcal{P})\subset E.\end{align*} Thanks to these properties, as
well as (P1),(P2) and (P3), we have the bound
\begin{multline}\label{8}|E|\gtrsim\iiint_{\mathcal{P}}
\left|J_{\Phi_{\mathbf{P}}}(s,t,u)\right|dsdtdu \gtrsim\\
\int_{S}|s-b|^{k/3}\int_{ T_{s}}|t-b|^{k/3}|s-t|
\int_{U_{s,t}}|u-b|^{k/3}|u-s||u-t|dudtds.\end{multline} To
estimate this integral from below, we shall have to split our
argument into three cases; our starting point will be to write
\[ T_{s}= T_{s}^1\cup T_{s}^2\cup T_{s}^3,\] where
\begin{align*}&  T_{s}^1=T_{s}\cap\{t\in I:|t-b|\leq (1/8)|s-b|\},\\
& T_{s}^2= T_{s}\cap\{t\in I:(1/8)|s-b|<|t-b|\leq 2|s-b|\},\\
& T_{s}^3= T_{s}\cap\{t\in I:|t-b|\geq 2|s-b|\}.\end{align*} Since
we are only guaranteed that one of the sets $ T_{s}^{\ell},\
\ell=1,2,3$ has $\mu$ measure at least $\beta$ (although two of
them or all of them might), it suffices to obtain the uniform
bound
\[\int_{S}\dots\int_{ T_{s}^{\ell}}\dots
\int _{U_{s,t}}\dots dudtds \ \ \gtrsim \ \ \alpha^4\beta^2,\]
under the assumption that $\mu(T_{s}^{\ell})\gtrsim \beta$ for
each $\ell = 1,2, 3$.\footnote{Strictly speaking, the choice of
$\ell\in \{1,2,3\}$ for which $T^{\ell}_s$ has large $\mu$ measure
depends on $s\in S$ and so, more accurately, we should split $S$
into three sets, stablising the choice of $\ell$ and noting that
one of these sets must have $\mu$ measure at least $\alpha$. We
hope our choice of exposition will not cause confusion.} However,
each case will be split into three subcases; to do so we shall
write
\[U_{s,t}=U_{s,t}^1\cup U_{s,t}^2\cup U_{s,t}^3,\]
  where
\begin{align*}& U_{s,t}^1=U_{s,t}\cap\{u\in I:|u-b|\leq
  (1/4)|t-b|\},\\
& U_{s,t}^2=U_{s,t}\cap\{u\in I:
  (1/4)|t-b|<|u-b|\leq 4|t-b|\},\\
& U_{s,t}^3=U_{s,t}\cap\{u\in I:|u-b|\geq
  4|t-b|\}.\end{align*}
Again, only one of the subsets $U_{s,t}^m,\ m=1,2,3$ is guaranteed
to have $\mu$ measure at least $\alpha,$ and our goal will be to
show the uniform bounds
\[\int_{S}\dots\int_{ T_{s}^{\ell}}\dots
\int _{U_{s,t}^m}\dots dudtds \ \ \gtrsim \ \ \alpha^4\beta^2,\]
under the assumptions $\mu(T_{s}^{\ell})\gtrsim \beta, \,
\mu(U_{s,t}^m)\gtrsim\alpha$ for each $\ell$ and $m =
1,2,3$.\footnote{Similar comments as above are valid here as
well.}

To successfully bound the iterated integral in \eqref{8} from
below we will need to excise various intervals from subsets of
$S$, $T^{\ell}_s$ and $U^m_{s,t}$ without changing their $\mu$
measure significantly. For this purpose we introduce the following
dynamic notation.
\begin{itemize}
\item For $\delta>0$, let $B_{\alpha} = \{u\in I: |u-b|\le \delta
\alpha^{6/(k+6)} \}$ so that $\mu(B_{\alpha}) \le c_k
\delta^{(k+6)/6} \alpha$. We will choose $\delta>0$ to be
sufficiently small in each instance so that the following holds:
if $W\subset I$ is a set satisfying $\mu(W)
> c_0 \alpha$ for some $c_0>0$, then $\mu(W\setminus B_{\alpha}) \ge (c_0/2) \alpha$
if $\delta>0$ is sufficiently small.

\item For $\delta>0$ and $t$, set $B_{t,\alpha} = \{u\in I:
|u-t|\le \delta \alpha |t-b|^{-k/6} \}$.
\begin{itemize}

\item If on $W\subset I$, $|u-b|\le C_0 |t-b|$, then $\mu(W\cap
B_{t,\alpha}) \le C_0^{k/6} \delta \alpha$ and therefore if
$\mu(W)\ge c_0 \alpha$, we have $\mu(W\setminus B_{t,\alpha}) \ge
(c_0/2) \alpha$ if $\delta>0$ is chosen sufficiently small.

\item On the other hand, if we do not know apriori that $|u-b| \le
C_0 |t-b|$ on $W$ but we happen to know $|t-b| \ge C_0
\alpha^{6/(k+6)}$, then automatically we have the control
$|u-b|\lesssim |t-b|$ on $B_{t,\alpha}$ since $|t-b| \ge C_0
\alpha^{6/(k+6)}$ implies $\alpha |t-b|^{-k/6} \lesssim |t-b|$ and
thus $|u-t| \lesssim |t-b|$ on $B_{t,\alpha}$.
\end{itemize}

\end{itemize}

\textbf{Case 1:} integration over $T_{s}^1;$ note that on this
set $|s-t|\sim |s-b|.$\\
Case 1a): integration over $U_{s,t}^1;$ here $|u-t|\sim |t-b|$ and
$|u-s| \sim |s-b|.$ Thus
\begin{multline*}\int_S |s-b|^{k/3}\int_{T^1_s} |t-b|^{k/3}|s-t|
\int_{U^1_{s,t}} |u-b|^{k/3}|u-s||u-t|dudtds \sim \\ \int_S
|s-b|^{k/3+2}\int_{T^1_s}|t-b|^{k/3+1}\int_{U^1_{s,t}} |u-b|^{k/3}dudtds\gtrsim\\
\int_{S\setminus B_{\alpha}} |s-b|^{k/6+k/6+2}\int_{T^1_s
\setminus B_{\beta}}|t-b|^{k/6+k/6+1} \int_{U^1_{s,t}\setminus
B_{\alpha}} |u-b|^{k/6+k/6}dudtds,
\end{multline*}
where we used the fact that on $U^1_{s,t} \setminus B_{\alpha}$ we
have $|u-b|\gtrsim \alpha^{\frac6{k+6}}$ (as well as analogous
estimates on $T^1_s \setminus B_{\beta}$ and $S\setminus
B_{\alpha}$). Now choosing $\delta>0$ in each $B_{\alpha},
B_{\beta}$ to ensure that the $\mu$ measure of the above sets have
not been altered significantly, we see that the last iterated
integral is bounded below by
\[\alpha^{\frac6{k+6}(k/6+2)}\times \alpha\times
\beta^{\frac6{k+6}(k/6+1)}\times\beta \times
\alpha^{\frac{k}{k+6}}\times \alpha =\alpha^4\beta^2.\]
Case 1b):
integration over $U_{s,t}^2;$ here $|u-s|\sim |s-b|$ but now
$|u-t|$ may vanish. Then
\begin{multline*}\int_S |s-b|^{k/3}\int_{T^1_s} |t-b|^{k/3}|s-t|
\int_{U^2_{s,t}} |u-b|^{k/3}|u-s||u-t|dudtds \sim \\ \int_S
|s-b|^{k/3+2}\int_{T^1_s}|t-b|^{k/3}\int_{U^2_{s,t}} |u-b|^{k/3}|u-t|dudtds\gtrsim\\
\int_S |s-b|^{k/3+2}\int_{T^1_s}
|t-b|^{k/3}\int_{U^2_{s,t}\setminus B_{t,\alpha}}
|u-b|^{k/3}|u-t|dudtds,\end{multline*} and using that on
$U^2_{s,t} \setminus B_{t,\alpha}$ one has $|u-t|\gtrsim\alpha
|t-b|^{-k/6}$ (together with the fact that $|s-b|\gtrsim |t-b|$
and $|u-b| \sim |t-b|$ in this case) this last quantity is bounded
below by
\begin{multline*}\alpha\int_S
|s-b|^{k/3+2}\int_{T^1_s} |t-b|^{k/6}\int_{U^2_{s,t}\setminus
B_{t,\alpha}}
|u-b|^{k/3}dudtds\gtrsim\\
\alpha\int_{S\setminus B_{\alpha}} |s-b|^{k/3+1}\int_{T^1_s
\setminus B_{\beta}}|t-b|^{k/3+1} \int_{(U^2_{s,t}\setminus
B_{t,\alpha})\setminus B_{\alpha}} |u-b|^{k/6}dudtds.
\end{multline*}
Since $|u-b|\le 2 |t-b|$ on $U^2_{s,t}$, we see that we can choose
$\delta>0$ in each $B_{\alpha}, B_{\beta}$ and $B_{t,\alpha}$ so
as not to change the $\mu$ measure much when we excise these
intervals from $S$, $T^1_s$ and $U^2_{s,t}$. Therefore the last
iterated integral above is at least
$\alpha\times\alpha^2\times\beta^2\times\alpha=\alpha^4\beta^2.$ \\
Case 1c): integration over $U_{s,t}^3;$ here $|u-t| \sim |u-b|$
but now $|u-s|$ may vanish. Then
\begin{multline*}\int_S |s-b|^{k/3}\int_{T^1_s} |t-b|^{k/3}|s-t|
\int_{U^3_{s,t}} |u-b|^{k/3}|u-s||u-t|dudtds\gtrsim\\
\int_{S\setminus B_{\alpha}} |s-b|^{k/3+1}\int_{T^1_s} |t-b|^{k/3}
\int_{U^3_{s,t}\setminus B_{s,\alpha}} |u-b|^{k/3+1}|u-s|dudtds\gtrsim\\
\alpha \int_{S\setminus B_{\alpha}} |s-b|^{k/6+1}\int_{T^1_s} |t-b|^{k/3}
\int_{U^3_{s,t} \setminus B_{s,\alpha}} |u-b|^{k/3+1}dudtds\gtrsim\\
\alpha \int_{S\setminus B_{\alpha}} |s-b|^{k/6}\int_{T^1_s
\setminus B_{\beta}} |t-b|^{k/3+1} \int_{(U^3_{s,t} \setminus
B_{s,\alpha})\setminus B_{\alpha}} |u-b|^{k/3+1}dudtds.
\end{multline*}
Since we do have the control $|u-b| \lesssim |s-b|$ on
$B_{s,\alpha}$ (since for $s\in S\setminus B_{\alpha}$, $|s-b|
\gtrsim \alpha^{6/(k+6)}$), we see that by appropriate choices of
$\delta>0$ in $B_{\alpha}, B_{\beta}$ and $B_{s,\alpha}$, the
above excised sets do not change in $\mu$ measure. Thus the final
iterated integral is at least
$\alpha\times\alpha\times\beta^2\times\alpha^2=\alpha^4\beta^2.$

\textbf{Case 2:} integration over $T_{s}^2.$\\
Case 2a): integration over $U_{s,t}^1;$ here $|u-t|\sim |t-b|,$
and we may also deduce $|u-s|\sim |s-b|.$ Since $|t-b| \sim
|s-b|$,
\begin{multline*}\int_S |s-b|^{k/3}\int_{T^2_s} |t-b|^{k/3}|s-t|
\int_{U^1_{s,t}} |u-b|^{k/3}|u-s||u-t|dudtds\gtrsim\\
\int_{S} |s-b|^{k/3+1}\int_{T^2_s \setminus
B_{s,\beta}} |t-b|^{k/3+1}|s-t| \int_{U^1_{s,t}}|u-b|^{k/3}dudtds \gtrsim\\
\beta \int_{S}
|s-b|^{k/6+1}\int_{T^2_s \setminus B_{s,\beta}} |t-b|^{k/3+1}
\int_{U^1_{s,t}} |u-b|^{k/3}dudtds
\gtrsim\\
\beta \int_{S\setminus B_{\alpha}} |s-b|^{k/3+2} \int_{T^2_s
\setminus B_{s,\beta}} |t-b|^{k/6} \int_{U^1_{s,t}\setminus
B_{\alpha}}|u-b|^{k/3}dudtds.
\end{multline*}
Again since $|t-b|\lesssim |s-b|$, appropriate choices of
$\delta>0$ can be made so as not to change the $\mu$ measure of
$S$, $T^2_s$ and $U^1_{s,t}$ when we excise from them the above
intervals. Hence the last iterated integral is at least
$$\beta\times\alpha\times \alpha^{\frac6{k+6}(k/6+2)}\times
\beta\times
\alpha\times\alpha^{\frac6{k+6}(k/6)}=\alpha^4\beta^2.$$\\
Case 2b): integration over  $U_{s,t}^2;$ here we may compare all
quantities containing $b$; namely $|s-b|\sim |t-b|\sim |u-b|.$
Hence
\begin{multline*}\int_S |s-b|^{k/3}\int_{T^2_s} |t-b|^{k/3}|s-t|
\int_{U^2_{s,t}} |u-b|^{k/3}|u-s||u-t|dudtds\gtrsim\\
\int_{S} |s-b|^{k/2}\int_{T^2_s \setminus
B_{s,\beta}}|t-b|^{k/3}|s-t| \int_{U^2_{s,t} \setminus
(B_{t,\alpha}\bigcup B_{s,\alpha})} |u-b|^{k/6} |u-t||u-s|
dudtds\gtrsim\\
\beta\alpha^2 \int_{S\setminus B_{\alpha}} |s-b|^{k/6}\int_{T^2_s
\setminus B_{s,\beta}}|t-b|^{k/6} \int_{U^2_{s,t} \setminus
(B_{t,\alpha}\bigcup B_{s,\alpha})} |u-b|^{k/6} dudtds.
\end{multline*}
Again we see that the sets we are integrating over have not
changed $\mu$ measure much when we remove intervals and so the
last iterated integral is at least
$\beta\alpha^2\times\alpha\times\beta\times
\alpha=\alpha^4\beta^2$.\\
Case 2c): integration over $U_{s,t}^3;$ here $|u-t|\sim |u-b|$ but
$|u-s|$ and $|t-s|$ may vanish. Since $|u-b| \gtrsim |s-b| \sim
|t-b|$,
\begin{multline*}\int_S |s-b|^{k/3}\int_{T^2_s} |t-b|^{k/3}|s-t|
\int_{U^3_{s,t}} |u-b|^{k/3}|u-s||u-t|dudtds\gtrsim\\
\int_{S\setminus B_{\alpha}} |s-b|^{k/3}\int_{T^2_s \setminus
B_{s,\beta}}  |t-b|^{k/3}|s-t|
\int_{U^3_{s,t} \setminus B_{s,\alpha}} |u-b|^{k/3+1} |u-s| dudtds\gtrsim\\
\alpha\beta \int_{S\setminus B_{\alpha}} |s-b|^{k/3 +
1}\int_{T^2_s \setminus B_{s,\beta}} |t-b|^{k/6} \int_{U^3_{s,t}
\setminus B_{s,\alpha}} |u-b|^{k/6}dudtds.
\end{multline*}
One checks that removing $B_{\alpha}$, $B_{s,\beta}$ and
$B_{s,\alpha}$ has not changed the $\mu$ measure of our sets very
much and so this last iterated integral is at least
$\alpha\beta\times\alpha^2\times\beta\times\alpha
=\alpha^4\beta^2$.

\textbf{Case 3:} integration over $T_{s}^3;$ in this interval
$|t-s| \sim |t-b|.$\\
Case 3a): integration over $U_{s,t}^1;$ here $|t-u| \sim |t-b|$
but $|u-s|$ may vanish. Since $|t-b|\gtrsim |s-b|$,
\begin{multline*}\int_S |s-b|^{k/3}\int_{T^3_s} |t-b|^{k/3}|s-t|
\int_{U^1_{s,t}} |u-b|^{k/3}|u-s||u-t|dudtds\gtrsim\\
\int_{S\setminus B_{\alpha}} |s-b|^{k/3}\int_{T^3_s \setminus
B_{\beta}} |t-b|^{k/3+2}\int_{U^1_{s,t} \setminus B_{s,\alpha}}
|u-b|^{k/3}|u-s| dudtds\gtrsim\\
\alpha \int_{S\setminus B_{\alpha}} |s-b|^{k/6+1}\int_{T^3_s
\setminus B_{\beta}} |t-b|^{k/3+1}\int_{(U^1_{s,t} \setminus
B_{s,\alpha}) \setminus B_{\alpha}} |u-b|^{k/3} dudtds.
\end{multline*}
Again the removal of intervals have not changed significantly the
$\mu$ measure and so the last iterated integral is at least $\alpha
\times\alpha^{k/(k+6) + 1}\times\beta^2\times\alpha^{6/(k+6) + 1}
= \alpha^4\beta^2$.\\
Case 3b): integration over $U_{s,t}^2;$ here $|s-u|\sim |u-b|$ but
$|u-t|$ can vanish. Since $|u-b| \sim |t-b|$,
 \begin{multline*}\int_S |s-b|^{k/3}\int_{T^3_s} |t-b|^{k/3}|s-t|
\int_{U^2_{s,t}} |u-b|^{k/3}|u-s||u-t|dudtds\gtrsim\\
\int_{S} |s-b|^{k/3}\int_{T^3_s} |t-b|^{k/3+1}
\int_{U^2_{s,t} \setminus B_{t,\alpha}} |u-b|^{k/3 +1}|u-t|dudtds\gtrsim\\
 \alpha \int_{S\setminus B_{\alpha}}|s-b|^{k/3}\int_{T^3_s
\setminus B_{\beta}} |t-b|^{k/3 +1} \int_{U^2_{s,t} \setminus
B_{t,\alpha}}|u-b|^{k/6 +1}dudtds,
\end{multline*}
and as before we see that the last iterated integral is at least
$\alpha \times \alpha^{k/(k+6) + 1} \times \beta^2 \times
\alpha^{6/(k+6) + 1} = \alpha^4 \beta^2.$\\
Case 3c): integration over $U_{s,t}^3;$ here we may deduce that
$|u-t| \sim |u-b|$ and $|u-s| \sim |u-b|.$ Thus
 \begin{multline*}\int_S |s-b|^{k/3}\int_{T^3_s} |t-b|^{k/3}|s-t|
\int_{U^3_{s,t}} |u-b|^{k/3}|u-s||u-t|dudtds\gtrsim\\
\int_{S\setminus B_{\alpha}} |s-b|^{k/3}\int_{T^3_s \setminus
B_{\beta}} |t-b|^{k/3+1} \int_{U^3_{s,t} \setminus B_{\alpha}}
|u-b|^{k/3+2}dudtds,
\end{multline*}
and as before this last iterated integral is at least
$\alpha\times\alpha^{\frac6{k+6}(k/6)}\times\beta^2\times\alpha\times
\alpha^{\frac6{k+6}(k/6+2)}=\alpha^4\beta^2.$

This completes the bound for \eqref{8} and thus proves
\eqref{rw-E}, completing the proof of Theorem
\ref{restricted}.\qed

\section{Strong-type inequalities}
We now wish to complete the proof of Theorem \ref{main} when
$d=3.$ We shall suitably modify the arguments in \cite{3} in order
to achieve this goal. We will concentrate only on the first
estimate stated in Theorem \ref{main} and thanks to our geometric
inequality and previous arguments, we just have to show that the
operator $\ca_I:L^2(\rthree)\to L^{3,2+\epsilon}(\rthree)$,
uniformly in $I$. This is equivalent to showing
\begin{equation}\label{9}
\left|\langle\ca_I f,g\rangle\right|\le C_{\epsilon}
\|f\|_2\|g\|_{3/2,2-\epsilon}\quad \text{any }f\in
L^2(\rthree),g\in L^{3/2,2-\epsilon}(\rthree).\end{equation}
Following \cite{3}, it suffices to select $f,g$ of the form
\[f=\sum_{\ell\in\zed}2^{\ell}\chi_{E_{\ell}},\quad
g=\sum_{m\in\zed}2^m \chi_{F_m},\] where the sets $E_{\ell}$ are
pairwise disjoint and so are the sets $F_m.$ However, we shall
specialise further, and pick the function $g = g_0$ to be simply
the characteristic function of a measurable set, $g_0 := \chi_F$.
If we prove estimate (\ref{9}) with $g$ replaced by $g_0,$ we then
have a $L^2\to L^{3,\infty}$ bound; one can then use Christ's
arguments to turn this into the claimed Lorentz space bound. We
may normalise the $L^2$ norm of $f,$ so that
$\sum_{\ell}2^{2\ell}|E_{\ell}|=1,$ and then the desired $L^2 \to
L^{3,\infty}$ bound becomes
\begin{equation}\label{weak-type}
\sum_{\ell} 2^{\ell} \langle\ca_I \chi_{E_{\ell}},\chi_F\rangle
\lesssim |F|^{2/3}.
\end{equation}
We decompose the $\ell$ sum above in order to stablise certain
quantities. For dyadic numbers $\epsilon,\eta\in(0,1/2]$ we define
$L_{\epsilon,\eta}$ to be those $\ell$ where
\[|E_{\ell}|\sim \eta 2^{-2\ell}\quad {\rm and} \quad
\langle\ca_I \chi_{E_{\ell}},\chi_F\rangle\sim \epsilon
|E_{\ell}|^{1/2} |F|^{2/3}.\] The number $M$ of indices $\ell$ in
$L_{\epsilon,\eta}$ is therefore finite and satisfies
$M\eta\lesssim 1.$ Our aim is then to prove
\begin{equation}\label{a,b}
\sum_{\ell\in L_{\epsilon,\eta}} 2^{\ell} \langle\ca_I
\chi_{E_{\ell}},\chi_F\rangle \lesssim \min(\epsilon^a, \eta^b)
|F|^{2/3}
\end{equation}
for some positive exponents $a,b$. By summing over the dyadic
$\epsilon$ and $\eta$, we see that \eqref{a,b} implies
\eqref{weak-type}.

Next we may assume that $|i-j|\geq C\log(1/\epsilon)$ for any two
distinct indices appearing in the sum over $L_{\epsilon,\eta}$
where $C>0$ will be an absolute constant. \footnote{By splitting
the sum over $L_{\epsilon,\eta}$ into $O(C\log(1/\epsilon))$ sums,
this assumption will cost us only a factor of
$O(C\log(1/\epsilon))$ in the estimate \eqref{a,b}.} One now
defines sets
\[G_{\ell}=\{x\in F:\ca_I\chi_{E_{\ell}}\geq c_0
|E_{\ell}|^{1/2}|F|^{2/3}|F|^{-1}\},\] for a certain $c_0>0$. If
$c_0$ is chosen sufficiently small, then $\langle\ca_I
\chi_{E_{\ell}},\chi_{F\setminus G_{\ell}} \rangle \le 1/2
\langle\ca_I \chi_{E_{\ell}},\chi_F\rangle$ and so $\langle\ca_I
\chi_{E_{\ell}},\chi_{G_{\ell}} \rangle \sim \langle\ca_I
\chi_{E_{\ell}},\chi_F\rangle$. By Theorem \ref{restricted} we
have $\langle\ca_I \chi_{E_{\ell}},\chi_{G_{\ell}} \rangle
\lesssim |E_{\ell}|^{1/2} |G_{\ell}|^{2/3}$ and so
\begin{equation}\label{G_k} |G_{\ell}|\gtrsim\epsilon^{3/2}|F|.
\end{equation} By the Cauchy-Schwarz inequality,
\begin{align*}
\bigl(|F|^{-1} \sum_{\ell\in L_{\epsilon,\eta}} |G_{\ell}|\bigr)^2
\le & |F|^{-1} \int_F \bigl(\sum_{\ell\in
L_{\epsilon,\eta}} \chi_{G_{\ell}}\bigr)^2 \\
\le & |F|^{-1} \sum_{\ell\in L_{\epsilon,\eta}}|G_{\ell}| +
|F|^{-1} \sum_{k\not= \ell} |G_k \cap G_{\ell}|
\end{align*}
and therefore either $\bigl(|F|^{-1} \sum_{\ell\in
L_{\epsilon,\eta}} |G_{\ell}|\bigr)^2
 \lesssim |F|^{-1} \sum_{k\not= \ell} |G_k \cap G_{\ell}|$ holds or we have \\
$\sum_{\ell\in L_{\epsilon,\eta}}|G_{\ell}|\lesssim
 |F|.$ If the former holds, then by \eqref{G_k}
$$
(M\epsilon^{3/2})^2 \lesssim \bigl(
\sum_{\ell\in
L_{\epsilon,\eta}} |G_{\ell}|\bigr)^2 \lesssim M^2 |F|^{-1}
\max_{k\not= \ell} |G_k \cap G_{\ell}|
$$
and the above dichotomy becomes
\begin{align}\label{10}&\text{either}\quad\sum_{\ell\in L_{\epsilon,\eta}}|G_{\ell}|\lesssim
 |F|\\
\label{11}&\text{or there exist }i\neq j\text{ so that} \quad
|G_i\cap G_j|\gtrsim \epsilon^3 |F|.\end{align} The key is now to
show that (\ref{11}) leads to a contradiction; this implies that
(\ref{10}) holds, and therefore
\begin{align*}
\sum_{\ell\in L_{\epsilon,\eta}} 2^{\ell} \langle\ca_I
\chi_{E_{\ell}},\chi_F\rangle \sim & \sum_{\ell\in
L_{\epsilon,\eta}} 2^{\ell} \langle\ca_I
\chi_{E_{\ell}},\chi_{G_{
\ell}} \rangle \\ & \lesssim \bigl(
\sum_{\ell \in L_{\epsilon,\eta}}
2^{3\ell}|E_{\ell}|^{3/2}\bigr)^{1/3} \bigl(\sum_{\ell\in
L_{\epsilon,\eta}} |G_{\ell}|\bigr)^{2/3} \lesssim \eta^{1/6}
|F|^{2/3}.
\end{align*}
On the other hand,
\begin{align*}
\sum_{\ell\in L_{\epsilon,\eta}} 2^{\ell} \langle\ca_I
\chi_{E_{\ell}},\chi_F\rangle \sim & \sum_{\ell\in
L_{\epsilon,\eta}} 2^{\ell} \epsilon |E_{\ell}|^{1/2}|F|^{2/3} \\
\lesssim & \ \ \epsilon M \eta^{1/2} |F|^{2/3} \lesssim \epsilon
\eta^{1/2} |F|^{2/3}
\end{align*}
and these two estimates together imply \eqref{a,b}.

To disprove (\ref{11}) we need the following result.
\begin{lem}\label{three-sets} Let $E,E',G\subset\rthree$ be measurable sets of finite
measure. Suppose that
\[\ca_I{\chi_E}(x)\geq\beta\quad\text{and}\quad
\ca_I{\chi_{E'}}(x)\geq\delta\quad\text{all } x\in G.\] If
$\beta'=\beta \frac{|G|}{|E|},$ then
\[|E'|\gtrsim \beta^A{\beta'}^2  \delta^B,\quad\text{with }1\leq A< 2,
\ 2< B\leq3, \ A+B=4.\]
\end{lem}
\noindent {\bf Proof } Set $\Phi_{\mathbf{P}}(s,t,u) = -
\mathbf{P}(s) + \mathbf{P}(t) - \mathbf{P}(u)$ and define
refinements
\begin{align*}&E^{1}=\left\{y\in E: \ca_I^*\chi_G(y)\geq \beta'/{2}\right\},\\
&G^{1}=\left\{x\in G: \ca_I\chi_{E^{1}}(x)\geq \beta/4\right\}.
\end{align*} We have
\begin{align*}
& \langle\ca_I^*\chi_{G^1},\chi_{E^1}\rangle =
\langle\ca_I\chi_{E^1},\chi_{G}\rangle - \langle\ca_I\chi_{E^1},\chi_{G\setminus G^1}\rangle  \ge  \langle\ca_I\chi_{E^1},\chi_{G}\rangle - \frac{\beta|G|}{4} \\
= \ &
\langle\ca_I^*\chi_{G},\chi_{E}\rangle-\langle\ca_I^*\chi_{G},\chi_{E\setminus
E^1}\rangle-\frac{\beta|G|}{4} \ge \langle\ca_I\chi_E ,\chi_G \rangle - \frac{3\beta|G|}{4} \ge \frac{\beta|G|}{4}.
\end{align*}
Hence $G^1\neq \emptyset.$ Now, pick $x_0\in G^1$ and
set
\[ S = \{s\in I:x_0-\mathbf{P}(s)\in E^1\}\Rightarrow \mu(S)=\ca_I\chi_{E^1}(x_0)\geq \beta/4.\]
For $s\in S$, set
\[T_{s}=\{t\in I:x_0-\mathbf{P}(s)+\mathbf{P}(t)\in G\}\Rightarrow \mu(T_{s})=\ca_I^*\chi_{G}(x_0-\mathbf{P}(s))\geq \frac{\beta|G|}{2|E|}.\]
Finally for $s\in S$ and $t\in T_s$, set
\[U_{s,t}=\{u\in I: x_0 + \Phi_{\mathbf{P}}(s,t,u) \in E'\} \Rightarrow \mu(U_{s,t})= \ca_I\chi_{E'}(x_0-\mathbf{P}(s)+\mathbf{P}(t))\geq \delta. \]
The idea is to estimate the measure of $E'$ by observing that if
\[\mathcal{P}=\left\{(s,t,u)\in I^3: s\in  S,t\in T_{s},
u\in U_{s,t}\right\}\text{ then } x_0+\Phi_{\mathbf{P}}
(\mathcal{P})\subset E'.\]
Hence the arguments of \S4 apply and we have
\[|E'|\gtrsim \int_{S}|s-b|^{k/3}\int_{T_s}|t-b|^{k/3}|s-t|
\int_{U_{s,t}}|u-b|^{k/3}|t-u||s-u| du dt ds,\] and this
quantity is bounded below by $\beta^A{\beta'}^2\delta^B,$ by the
proof of Theorem \ref{restricted}; note that the relation between
the numbers $A$ and $B$ can also be easily extracted from there. \qed

We can now conclude our argument; pick $E=E_i, E'=E_j,G=G_i\cap
G_j,$ and $\beta=\epsilon |E|^{1/2}|F|^{-1/3}, \delta=
|E'|^{1/2}|F|^{-1/3}, \beta'=\beta |G|/{|E|}.$ By Lemma
\ref{three-sets} we have
\begin{multline*}|E'|\gtrsim
\epsilon^{A+B}|F|^{(A+B)/3}|E|^{A/2}|F|^{B/2}\beta^2|G|^2|E|^{-2}
\gtrsim\\ \epsilon^4|F|^{-4/3}|E|^{A/2}|E'|^{B/2}\epsilon^2|E||F|^{-2/3}
|G|^2|E|^{-2}\gtrsim \epsilon^{12} |E|^{A/2-1}|E'|^{B/2},\end{multline*}
where we have used the fact that $|G|\gtrsim \epsilon^{3}|F|.$
Using the relation $A+B=4$ we deduce
\[|E'|^{1-A/2}\lesssim\epsilon^{-12}|E|^{1-A/2}\iff 2^{-jp}\lesssim
\epsilon^{24/{(2-A)}}2^{-ip}\] implying $j\geq
i-C'\log(1/{\epsilon});$ since the roles of $i,j$ can be exchanged
one has $|i-j|\leq C' \log(1/{\epsilon}),$ which contradicts our
assumptions and therefore \eqref{11} cannot hold. This gives us
the weak-type bound \eqref{weak-type}. As we have already
mentioned, the arguments in \cite{3} can now be reproduced
verbatim to obtain the Lorentz bound (\ref{9}), completing the
proof of Theorem \ref{main} for $d=3.$

\section{Two-dimensional estimates}
In this section we present the arguments necessary to prove
Theorem \ref{main} in the case $d=2,$ starting with the restricted
weak type estimates.
\begin{thm}\label{rw-2} Let $d=2.$ The operator (\ref{5}) satisfies
\begin{equation}\ca_\R :L^{3/2,1}(\rtwo)\to L^{3,\infty}(\rtwo).\end{equation}
\end{thm}
\noindent {\bf Proof }The preparatory statements of \S3 and \S4 can obviously be applied
also in this setting and we quickly reduce our analysis to the
operators
\[\ca_{I}f(x)=\int_I f(x-\mathbf{P}(t))|t-b|^{k/3}dt := \int_I
f(x-\mathbf{P}(t))d\mu_I(t),\]for each fixed $I.$ We set
\[\langle\ca_I \chi_E,\chi_F\rangle=\alpha_I |F|,\quad
\langle\ca_I \chi_E,\chi_F\rangle=\beta_I |E|,\] with $|E|\neq 0,
|F|\neq 0,$ and observe it suffices to establish\footnote{We shall again abuse notation and relabel the measure
$\mu_I$ as $\mu;$ the numbers $\alpha_I,\beta_I$ will also be
replaced by $\alpha$ and $\beta$.}
\begin{equation}\label{main-estimate}
\langle\ca_I \chi_E,\chi_F\rangle\lesssim |E|^{2/3}|F|^{2/3}\iff
|E|\gtrsim \alpha^2\beta\iff |F|\gtrsim\beta^2\alpha,
\end{equation}
uniformly in $I$. As discussed in \S2 we will apply Christ's
argument to prove
\begin{equation}\label{main-again}
|E|\gtrsim \alpha^2\beta \ {\rm only \ for} \ \alpha\le\beta \ \
{\rm and \ similarly} \ \ |F|\gtrsim\beta^2\alpha \ {\rm only \
for} \ \beta\le \alpha.
\end{equation}
But from the relation $\alpha |F| = \beta |E|$, we see that
\eqref{main-again} implies \eqref{main-estimate}. This only works
since we are proving an estimate on the line of duality. We shall
concentrate on the first estimate in \eqref{main-again} (the proof
of the second estimate is similar) and so we assume from now on
that $\alpha \le \beta$.

By the discussion in \S2 we can find a point $x_0\in E$ and
\begin{align*}& S\subset I\text{ so that
  }\mu(S)\gtrsim\beta; \\
& \text{for each }s\in S \text{ there is }T_{s}\subset
  I\text{ so that }\mu(T_{s})\gtrsim\alpha;\\
&\text{if }\mathcal{P}=\left\{(s,t)\in I^2: s\in  S,t\in T_{s},
\right\}\Longrightarrow x_0+\Phi_{\mathbf{P}} (\mathcal{P})\subset
E.\end{align*} Therefore (see \S2)
\begin{equation}\label{13}|E|\gtrsim\iint_{\mathcal{P}}
\left|J_{\Phi_{\mathbf{P}}}(s,t)\right|dsdt \gtrsim
\int_{S}|s-b|^{k/2}\int_{ T_{s}}|t-b|^{k/2}|s-t|
dtds.\end{equation} We split
\[ T_{s}= T_{s}^1\cup T_{s}^2\cup T_{s}^3,\]
where
\begin{align*}&  T_{s}^1=T_{s}\cap\{t\in I:|t-b|\leq (1/2)|s-b|\},\\
& T_{s}^2= T_{s}\cap\{t\in I:(1/2)|s-b|<|t-b|\leq 2|s-b|\},\\
& T_{s}^3= T_{s}\cap\{t\in I:|t-b|\geq 2|s-b|\}.\end{align*} By
the same arguments of \S4 we shall prove the bound
\[ \int_{S}|s-b|^{k/2}\int_{ T_{s}^{\ell}}|t-b|^{k/2}|s-t|
dtds\gtrsim\alpha^2\beta,\quad\ell=1,2,3\] under the assumption
that $\mu(T_{s}^\ell)\gtrsim\alpha$ in each case.

We shall use similar dynamic notation as in \S4: $B_{\alpha} =
\{t\in I: |t-b|\le \delta \alpha^{3/(k+3)} \}$ and $B_{s,\alpha} =
\{t\in I: |t-s| \le \delta\alpha |s-b|^{-k/3} \}$ with analogous
conclusions as before if $\delta>0$ is chosen small enough in any
particular situation.

\textbf{Case 1:} integration over $T_{s}^1;$ here $2|t-b| \le |s-b| \sim |t-s|.$ Thus
\[
\int_S |s-b|^{k/2}\int_{T^1_s} |t-b|^{k/2}|s-t| dtds\gtrsim
\int_{S \setminus B_{\beta}}|s-b|^{k/2+1}\int_{T^1_s \setminus
B_{\alpha}} |t-b|^{k/2} dtds \gtrsim \]
\[
\int_{S \setminus B_{\beta}}|s-b|^{k/3}\int_{T^1_s \setminus
B_{\alpha}} |t-b|^{2k/3+1} dtds \gtrsim \beta \alpha^2.
\]
Here we have not used the relation $\alpha \le \beta$. In addition,
\[ \int_{S \setminus B_{\beta}}|s-b|^{k/2+1}\int_{T^1_s \setminus
B_{\alpha}} |t-b|^{k/2} dtds \gtrsim \beta^{\frac3{2}\frac{k+4}{k+3}}
\alpha^{\frac{3}{2}\frac{k+2}{k+3}}. \]
Notice that $\beta^{\frac3{2}\frac{k+4}{k+3}}
\alpha^{\frac{3}{2}\frac{k+2}{k+3}}\gtrsim\alpha^2\beta$ for
$\alpha \le \beta$.
The former of these two estimates suffices for the proof of Theorem \ref{rw-2}. However, both estimates will be required in order to obtain Lorentz space bounds.

\textbf{Case 2:} integration over $T_{s}^2;$ we have

$$
\int_S |s-b|^{k/2}\int_{T^2_s} |t-b|^{k/2}|s-t| dtds\gtrsim
\int_{S\setminus B_{\beta}} |s-b|^{k/2}\int_{T^2_s \setminus
B_{s,\alpha}} |t-b|^{k/2}|s-t| dtds\gtrsim
$$
$$
\alpha \int_{S\setminus B_{\beta}} |s-b|^{k/6}\int_{T^2_s
\setminus B_{s,\alpha}}|t-b|^{k/2} dtds.
$$
We make the important observation here that, in this case, $|t-b|\lesssim |s-b|$ on $B_{s,\alpha}$  and
therefore $\mu(T^2_s \setminus B_{s,\alpha}) \gtrsim \alpha$ if
$\delta > 0$ is chosen appropriately. Therefore the last iterated
integral is bounded below by
$$
\alpha \int_{S\setminus B_{\beta}}|s-b|^{k/6+k/6}\int_{T^2_s
\setminus B_{s,\alpha}}|t-b|^{k/3} dtds\gtrsim \alpha^2\beta.
$$

\textbf{Case 3:} integration over $T_{s}^3;$ here $|t-s|\sim
|t-b|.$ Thus
$$
\int_S |s-b|^{k/2}\int_{T^3_s}|t-b|^{k/2}|s-t| dtds\gtrsim \int_S
|s-b|^{k/2}\int_{T^3_s} |t-b|^{k/2+1}dtds \gtrsim
$$
$$
\int_{S\setminus B_{\beta}} |s-b|^{k/2}\int_{T^3_s \setminus
B_{\alpha}} |t-b|^{k/2+1}dtds \gtrsim
\beta^{\frac{3}{2}\frac{k+2}{k+3}}
\alpha^{\frac{3}{2}\frac{k+4}{k+3}}\gtrsim \alpha^2\beta
$$
since $\alpha\le \beta$. This completes the proof of
\eqref{main-again} and hence the proof of Theorem \ref{rw-2}. \qed

To prove the Lorentz estimates for the operator $\ca_I$ we put
ourselves back in the setting of \S5, with the (obvious)
difference that we must consider the estimates just proven. Recall
the appropriate setup:
\begin{itemize}\item[-] there are 4 sets $E (= E_i), E' (= E_j), G (= G_i \cap G_j) , F$ with
$|E|\sim \eta2^{-3i/2},\ |E'|\sim \eta 2^{-3j/2},$ and
$G\subset F,$
\item[-] four parameters $\epsilon>0,\ \beta=\epsilon|E|^{2/3}|F|^{-1/3}
,\ \delta=\epsilon|E'|^{2/3}|F|^{-1/3},\ \beta'=\beta{|G|}/{|E|},$
\item[-] we may assume $|G|>\epsilon^3 |F|,\ \ca_I\chi_{E}\gtrsim\beta$
  on $G,$ $\ca_I\chi_{E'}\gtrsim\delta$ on $G,$
\item[-] we further assume $\beta\leq\delta\iff |E|\leq
|E'|,$\footnote{Since our arguments are completely symmetrical,
this assumption does not pose any restrictions, as the roles of
$E$ and $E'$ can be interchanged.}
\end{itemize}
and we wish to show that (\ref{11}) leads to a contradiction;
this will manifest itself in two possible forms, the inequality
\[|E|\gtrsim \epsilon^c |E'|\quad\text{or the inequality}\quad
|G|\le K^{-1}\epsilon^3  |F|,\] for some $c\ge 0$ and for a
sufficiently large $K$. Clearly $|G|\le K^{-1}\epsilon^3  |F|$
contradicts (\ref{11}). The inequality $|E|\gtrsim \epsilon^c
|E'|$ is equivalent to $2^{3(i-j)/2} \lesssim (1/\epsilon)^c$
which in turn is equivalent to $0 \le i-j \lesssim c\log
(1/\epsilon)$ which contradicts our basic assumptions on $i$ and
$j$. As indicated at the end of \S2 the arguments in \S5 break
down in the two dimensional setting and a slightly more elaborate
argument is needed here. To carry out our arguments, we define two
refinements
\[E^{1}=\left\{y\in E: \ca_I^*\chi_G(y)\geq \beta'/2\right\},\quad
 G^{1}=\left\{x\in G: \ca_I\chi_{E^{1}}(x)\geq \beta/4\right\}.
\]
The standard argument shows that $G^1\neq\emptyset,$ thus we pick
$x_0\in  G^1$ and set
\[S = \{s\in I:x_0-\mathbf{P}(s)\in
  E^1\}\Rightarrow\mu(S)=\ca_I\chi_{E^1}(x_0)\geq\beta/{4},\]
\[ T_{s}=\{t\in I:x_0-\mathbf{P}(s)+\mathbf{P}(t)\in
  G\}\Rightarrow \mu(T_{s})=\ca_I^*\chi_{G}(x_0-\mathbf{P}(s))\geq
 {\beta'}/{2},\]
\begin{multline*}
 U_{s,t}=\{u\in I:x_0 -\mathbf{P}(s)+\mathbf{P}(t)-\mathbf{P}(u) \in E'\} \\
 \Rightarrow \mu(U_{s,t})=\ca_I\chi_{E'}(x_0-\mathbf{P}(s)+\mathbf{P}(t))\geq
\delta. \end{multline*}

\textbf{Case A:} $|G|\geq\epsilon^p |E|,$ where $p>0$ will be
determined later.\\
For fixed $s\in S$ we have
\[ \psi_s(T_{s}\times U_{s,t}) \subset E',\quad\text{where}\quad
\psi_s(t,u)=x_0-\mathbf{P}(s)+\mathbf{P}(t)-\mathbf{P}(u),\] therefore
\[|E'|\gtrsim\int_{T_{s}}|t-b|^{k/2}\int_{U_{s,t}}|u-b|^{k/2}|u-t|dudt
\gtrsim \delta^C{\beta'}^D\] thanks to Cases 1, 2 and 3 in this
section; here $(C,D)=(2,1)$, $(A, B)$ or $(B, A)$, where
$(A, B) :=(\frac{3}{2}\frac{k+4}{k+3},\frac{3}{2}\frac{k+2}{k+3})$, and in all instances $C+D=3.$ Hence
\begin{multline*}|E'|\gtrsim \delta^C{\beta'}^D=
\delta^C\beta^D|G|^D|E|^{-D}\geq
\delta^C\beta^D\epsilon^{p(D-1)}|G||E|^{-1}\gtrsim\\
\epsilon^C|E'|^{2C/3}|F|^{-C/3}\epsilon^D|E|^{2D/3}|F|^{-D/3}
\epsilon^{p(D-1)}|G||E|^{-1},\end{multline*}
which is equivalent to
\[|E|^{1-2D/3}\gtrsim\epsilon^{3+p(D-1)}|E'|^{2C/3-1}|F|^{-1}|G|\gtrsim
\epsilon^{6+p(D-1)}|E'|^{2C/3-1},\]
the contradiction we wished to find.

\textbf{Case B:} $|G|\leq\epsilon^p|E|.$ This case is more
involved and will be split into subcases. However, we shall not
change our setup. Let
\[\mathcal{Q}=\{(s,t)\in I^2: s\in S, \ t\in T_{s}\},\quad
\Phi_{\mathbf{P}}(s,t)=x_0-\mathbf{P}(s)+\mathbf{P}(t).\]
Clearly $\Phi_{\mathbf{P}}(\mathcal{Q})\subset G,$ hence
\[|G|\gtrsim \int|s-b|^{k/2}\int|t-b|^{k/2}|s-t|dtds.\] Let
\[T_{s}=T_{s}^1\cup T_{s}^2\cup T_{s}^3,\]
where the sets $T_{s}^{\ell}$, $\ell=1,2,3$ are defined as above.
Also let
\[S^1=\{s\in S:\mu(T^2_{s})\geq\beta'/{6}\}, \quad
S^2=\{s\in S:\mu(T^1_{s})\geq\beta'/{6}\},\]
\[S^3=\{s\in S:\mu(T^3_{s})\geq\beta'/{6}\}.\]

\textsf{Case B1:} $\mu(S^1)\leq\beta/{12}.$ Then either
$\mu(S^2)\geq\beta/{12}$ or
$\mu(S^3)\geq\beta/{12}.$\\
Case B1a):  $\mu(S^2)\geq\beta/{12}.$ In this case, by Case 1,
\[
|G|\gtrsim\int_{S^2}|s-b|^{k/2}\int_{T^1_{s}}|t-b|^{k/2}|s-t|dtds
\gtrsim \beta^A{\beta'}^B=\epsilon^3 |E|^2|F|^{-1}(|G|/{|E|})^B.
\]
This implies
\[|F|\gtrsim \epsilon^3|E|^{2-B}|G|^{B-1}\geq \epsilon^3
\epsilon^{-p(2-B)}|G|^{2-B+B-1}\iff
|G|\lesssim\epsilon^{p(2-B)-3}|F|,\]
contradicting $|G|\gtrsim\epsilon^3|F|$ for $p$ chosen sufficiently large (note $B<2$).\\
Case B1b):  $\mu(S^3)\geq\beta/{12}.$ Here, by Case 3,
\[
|G|\gtrsim\int_{S^3}|s-b|^{k/2}\int_{T^3_{s}}|t-b|^{k/2}|s-t|dtds
\gtrsim {\beta'}^A{\beta}^B=\epsilon^3
|E|^2|F|^{-1}(|G|/{|E|})^A.
\]
This leads to
\[|F|\gtrsim \epsilon^3|E|^{2-A}|G|^{A-1}\geq \epsilon^3
\epsilon^{-p(2-A)}|G|^{2-A+A-1}\iff
|G|\lesssim\epsilon^{p(2-A)-3}|F|,\]
contradicting $|G|\gtrsim\epsilon^3|F|$ for sufficiently large $p$ (note $A<2$ if $k\not=0$\footnote{The case $k=0$ is simpler and is dealt with in \cite{3}.}).

\textsf{Case B2:}  $\mu(S^1)>\beta/{12}.$ To take care of this
case we shall define subsets $T^{2,1}_{s},T^{2,2}_{s}$
 of $T^2_{s}$ as
\begin{gather*}T^{2,1}_{s}=\{t\in T^{2}_{s}:\mu(\{u\in
U_{s,t}:|u-b|\leq 2|t-b|\})\geq\delta/2\},\\
T^{2,2}_{s}=\{t\in T^{2}_{s}:\mu(\{u\in U_{s,t}:|u-b|>
2|t-b|\})\geq\delta/2\}.
\end{gather*}
Case B2a): there exists $s_0\in S^1$ so that
$\mu(T^{2,1}_{s_0})\geq \beta'/{12}.$ Hence, we bound the measure
of $E'$ by integrating over $T^{2,1}_{s_0}$. By Cases 1 and 2, we have
\[|E'|\gtrsim \int_{T^{2,1}_{s_0}} |t-b|^{k/2}\int_{U_{s_0 ,t}} |u-b|^{k/2}|u-t|
  dudt
\gtrsim
\beta'\delta^2=\epsilon^3|E|^{-1/3}|E'|^{4/3}|G||F|^{-1}.\]
This implies $|E|^{1/3}\gtrsim \epsilon^3 |E'|^{1/3}|G||F|^{-1}\gtrsim \epsilon^6
|E'|^{1/3}$, giving us the desired contradiction. \\
Case B2b): for every $s\in S^1$ we have
$\mu(T^{2,1}_{s})< \beta'/{12}.$ Thus, we must have that
$\mu(T^{2,2}_{s})\geq \beta'/{12}.$ Now the integration occurs
over $T^{2,2}_{s};$ fixing an $s\in S^1$, we have
\[
|E'|\gtrsim |s-b|^{k/6}\int_{T^{2,2}_{s}} |t-b|^{k/3}\int_{U_{s, t}} |u-b|^{k/2}|u-t|dudt \gtrsim |s-b|^{k/6}\delta^A \beta'.
\]
Now, if we choose $\mathfrak{S}\subset S^1,$ so that
$\mu(\mathfrak{S})=\beta/{100}$ we have
\[|E'|\int_{\mathfrak{S}}|s-b|^{k/3}ds\gtrsim
\delta^A\beta'\int_{\mathfrak {S}\setminus B_{\beta}}
|s-b|^{k/3+k/6}ds\gtrsim\delta^A\beta'\beta\times
\beta^{\frac{k}{6}\frac3{k+3}}=\delta^A\beta'\beta^B,
\] and this implies
\begin{align*}
& \beta|E'|\gtrsim \delta^A\beta'\beta^B\iff |E'|\gtrsim\delta^A\beta^B
|G||E|^{-1}=\epsilon^3|E'|^{2A/3}|E|^{2B/3-1}|F|^{-1}|G| \\
\geq \ & \epsilon^6|E'|^{2A/3}|E|^{2B/3-1}\iff |E|^{1-2B/3}\gtrsim\epsilon^6
|E'|^{2A/3-1},
\end{align*}
which is the required contradiction.

\section{Sharpness of Theorem \ref{main}}
In this section we wish to show how the result of Theorem
\ref{main} is essentially sharp in the scale of Lorentz spaces by
providing an explicit, possibly well-known counterexample.
Consider the translation invariant operator $S$ given by
\[Sf(x)=\int_{-1}^1 f(x_1-t,x_2-t^2,\ldots, x_d-t^d) dt,\]
along with the family of nonisotropic dilations
\[\delta\circ y=(\delta y_1,\delta^2 y_2,\ldots,\delta^d y_d),\quad
\delta>0,\ y\in\rd.\]
 For a positive integer $k$, we let $K=(k,k^2,\ldots,k^d)$ and $\chi=\chi(y)=
\chi_{[-1/2,1/2]^d}(y).$ Further, let $\chi_k=\chi_k(y)\equiv\chi(k\circ y).$
For $N$ chosen sufficiently large, we define
\[f(x)=\sum_{k\geq N}\chi_k(x-K).\]
The supports $E_k$ of the functions appearing in the sum are disjoint for large enough $N$. Hence
\[\|f\|_{\frac{d+1}{2}}=\left(\sum_{k\geq N}
  |E_k|\right)^{\frac{2}{d+1}}\sim
\left(\sum_{k\geq N}k^{\frac{-d(d+1)}{2}}\right)^{\frac{2}{d+1}}\sim
N^{(1-\frac{d(d+1)}{2})\frac{2}{(d+1)}}=N^{\frac{2}{d+1}-d}.\]
Now
\begin{multline*}Sf(x)=\sum_{k\geq
  N}\int_{-1}^1\chi_k(x_1-t-k,\ldots,x_d-t^d-k^d)dt\geq\\
\sum_{k\geq
  N}\int_{|t|\leq k^{-1}/{10}}\chi_k(x_1-t-k,\ldots,x_d-t^d-k^d)dt\gtrsim
\sum_{k\geq N} k^{-1}\chi_{2k}(x-K),\end{multline*}
where the functions involved  in the last sum again have disjoint
supports $F_k.$ Thus, we may deduce
\begin{align*}& \|Sf\|_{L^{\frac{d+1}{2}\frac{d}{d-1},r}}\gtrsim
\left[\sum_{k\geq
    N}\left(k^{-1}|F_k|^\frac{2(d-1)}{d(d+1)}\right)^r\right]^{1/r}=
\left[\sum_{k\geq N}\left(k^{-1}
k^{-\frac{d(d+1)}{2}\frac{2(d-1)}{d(d+1)}}\right)^r\right]^{1/r} \\
& = \left[\sum_{k\geq N}k^{-dr}\right]^{1/r}\sim
(N^{-dr+1})^{1/r}= N^{-d+1/r}.
\end{align*} Hence, in order to
have boundedness, we must have the inequality
\[N^{-d+1/r}\lesssim N^{2/{(d+1)}-d},\]
which for sufficiently large $N$ implies
\[-d+1/r\leq 2/{(d+1)}-d\iff
    r\geq (d+1)/{2}.\]
Hence the result of Theorem \ref{main} is indeed sharp, except
possibly for the appearance of the $\epsilon$.

\end{document}